\documentclass[12pt]{article}
\usepackage{amsmath,amssymb,amsthm,latexsym}

\newtheorem{df}{Definition}
\newtheorem{lemma}{Lemma}

\newtheorem{theorem}{Theorem}
\newtheorem{prop}{Proposition}

\numberwithin{equation}{section}

\begin{document}
\title{{\sc Wrapping Brownian motion and heat kernels  I: compact Lie groups}}

\author{David G. Maher}
%\date{}
\maketitle

\newcommand{\g}{\mathfrak{g}}
\newcommand{\kg}{\mathfrak{k}}
\newcommand{\tg}{\mathfrak{t}}
\newcommand{\ug}{\mathfrak{u}}
\newcommand{\p}{\mathfrak{p}}
\newcommand{\X}{\mathfrak{X}}
\newcommand{\af}{\mathfrak{a}}
\newcommand{\R}{\mathbb{R}}
\newcommand{\C}{\mathbb{C}}
\newcommand{\Z}{\mathbb{Z}}
\newcommand{\N}{\mathbb{N}}
\newcommand{\E}{\mathbb{E}}
\newcommand{\bbP}{\mathbb{P}}
\newcommand{\T}{\mathbb{T}}
\newcommand{\F}{\mathcal{F}}
\newcommand{\A}{\mathcal{A}}
\newcommand{\x}{\mathbf{x}}
\newcommand{\y}{\mathbf{y}}
\newcommand{\ad}{\mathrm{ad}}
\newcommand{\Ad}{\mathrm{Ad}}
\newcommand{\Exp}{\mathrm{Exp}}
\newcommand{\grad}{\mathrm{grad} \,}

\begin{abstract}

An important object of study in harmonic analysis is the heat
equation. On a Euclidean space, the fundamental solution of the
associated semigroup is known as the heat kernel, which is also the
law of Brownian motion. Similar statements also hold in the case of
a Lie group. By using the wrapping map of Dooley and Wildberger, we
show how to wrap a Brownian motion to a compact Lie group from its
Lie algebra (viewed as a Euclidean space) and find the heat kernel.
This is achieved by considering It\^o type stochastic differential
equations and applying the Feynman-Ka\v{c} theorem.

\end{abstract}

{\small {\it Keywords:} Compact Lie group, Brownian motion, wrapping
map.

}

\section{Introduction}

The partial differential equation given on $\R^n$ by

\begin{equation}\label{311}
\partial_t u(x, t) = \tfrac{1}{2} \Delta u(x, t), \phantom{abcde} t \in \R^+, \; x \in \R^n,
\end{equation}
where $\Delta$ is the Laplacian, represents the dissipation of heat
over a certain time.  The fundamental solution of the associated
semigroup $e^{t \Delta / 2}$, known as the heat kernel, $p_t$ is
given by a unique, strongly continuous, contraction semigroup of
convolution operators which may be convolved with the initial data
$f(x) = u(0,x)$ to give the solution to the Cauchy problem.  The
heat kernel may also be expressed as the transition density of a
Brownian motion, $B_t$.  In summary:
\begin{align}
u(x,t) & = e^{t \Delta / 2} f(x)\\
& = (p_t * f)(x)\\
& = \int_{\R^n} p_t(x-y) f(y) dy\\
& = \E(f(B_t))
\end{align}

Similar statements hold when $\R^n$ is replaced by a Lie group.\\

In this article we will demonstrate how these results may be
transferred from the Lie algebra (regarded as $\R^n$) to a compact
Lie group using the so-called wrapping map of Dooley and Wildberger.
Additionally, we shall provide the mechanism that allows one to
``wrap" a Brownian motion, and then find the heat kernel by taking
the expectation of the ``wrapped" process and applying a
Feynman-Ka\v{c} type transform.\\

Results concerning Brownian motion and heat kernels a compact Lie
group have been previously given by many authors (eg, \cite{ARE},
\cite{ESK}, \cite{URA} and \cite{WAT} to name but a few).  Our
method is quite different by using the wrapping map, which can be
viewed as a global version of the exponential map.   As a result,
our results are obtained are in the spirit of the tangent space
analysis advocated by Helgason (\cite{HE2}, \cite{HE3}).\\

Finally, we will also discuss how our results may be extended and
particularly to compact symmetric spaces and complex Lie groups,
whose results will comprise of the sequel.

\section{Acknowledgements}

These results were obtained during the authors Ph.D candidature at
the University of New South Wales (see \cite{M}).  I would like to
thank my supervisor, Tony Dooley, for all his support and guidance
during this time.

\section{Notation}

Let $G$ be a compact, connected, simply connected Lie group, $\g$
its Lie algebra, $T$ a maximal torus and $\tg$ the Lie algebra of
$T$.  Let $n$ be the dimension of $G$, and $l$ the dimension of $T$
(also known as the rank of $G$).\\

Let $\g^\C = \g \otimes \C$ be the complexification of $\g$, $\tg^\C
= \tg \otimes \C$ be the complexification of $\tg$, and $\tg^\C$ is
a Cartan subalgebra of $\g^\C$.  Let $B( \cdot , \cdot )$ be the
Killing form on $\g^\C \times \g^\C$, with $\g^*$ and $\tg^*$ the
respective duals of $\g$ and $\tg$ with respect to the Killing
form.\\

We denote by $\Sigma$ the set of roots of $(\g^\C,\tg^\C)$, and
choose an ordering on $\Sigma$, with $\Sigma^+$ the set of positive
roots, $W$ the Weyl group, $\tg^*_+$ the positive Weyl chamber, and
let $\rho = \sum_{\alpha \in \Sigma^+} \alpha$. Let $k = \bigl|
\Sigma^+ \bigr|$ be the number of positive; thus we have $n = l +
2k$.  We denote by $\Lambda \subseteq \tg^*$, and let $\Lambda^+ =
\Lambda \cap \tg^*_+$ denote the set of positive integral weights.\\

We normalise Haar measures $dg$ on $G$ and $dt$ on $T$ to have total
mass 1.  Lebesgue measure $dX$ on $\g$ is normalised so that if $U$
is a neighbourhood of $0 \in \g$ on which the exponential map is
injective, then for $f \in C^\infty (G)$:
$$
\int_{U} f(\exp X) |j(X)|^2 dX = \int_{\exp U} f(g) dg:
$$
where $j(X)$ is is the analytic square root of the Jacobian or
$\exp$ with $j(0) = 1$, given by
$$
j(X) = \prod_{\alpha \in \Sigma^+} \frac{\sin \alpha (X)/2}{\alpha
(X)/2}
$$

Every irreducible representation $\pi \in \widehat{G}$ is associated
to a unique highest weight $\lambda \in \lambda^+$.  If $\chi_\pi =
\chi_\lambda$ is the character of this representation, the
Kirillov's character formula is given by
$$
j(X) \chi_\lambda (\exp X) \int_{\mathcal{O}_{\lambda + \rho}}
e^{i\beta(X)} d\mu_{\lambda + \rho} (\beta), \; \; \; \text{for all}
\; \; \; X \in \g
$$
where $\mathcal{O}_{\lambda + \rho}$ is the co-adjoint orbit through
$\lambda + \rho \in \tg^*_+$, and $\mu_{\lambda + \rho}$ is the
Liouville measure on $\mathcal{O}_{\lambda + \rho}$ with total mass
$d_\lambda = d_\pi = \dim \pi$.\\

\section{The wrapping map}

The {\it wrapping map}, $\Phi$, was devised by Dooley and Wildberger
in \cite{DW2}.  $\Phi$ is defined by
\begin{equation}
\langle \Phi(\nu), f \rangle = \langle \nu,j \tilde{f} \rangle
\end{equation}
where $f \in C^\infty (G)$, $\tilde{f} = f \circ \exp$ and $j$ the
analytic square root of the determinant of the exponential map.  We
need to place some conditions on $\nu$ for $\Phi(\nu)$ to be
well-defined  -  this is the case when $\nu$ is a distribution of
compact support on $\g$, or $j\nu \in L^1(\g)$ .  We call
$\Phi(\nu)$ the \textit{wrap} of $\nu$. When $\varphi$ is an
$G$-invariant Schwartz function, we have the following:
\begin{theorem}\label{thm1} (\cite{DW2}, Thm. 1) Let $\varphi \in S(\g)$ be $G$-invariant, then $\Phi(j \varphi)$ is a $G$-invariant $C^\infty$ function on $G$ given on $T$ by
\begin{equation}
\Phi(j \varphi) (\exp \, H) = \sum_{\gamma \in \Gamma} \varphi (H +
\gamma) \; \; \; \forall H \in \tg.
\end{equation}
\end{theorem}

The principal result is the {\bf wrapping formula}, given by
\begin{theorem}\label{thm2} (\cite{DW2}, Thm. 2) {Let $\mu, \nu$ be $G$-invariant distributions of compact support on $\g$ or two $G$-invariant integrable functions, then
\begin{equation}\label{wrapform}
\Phi(\mu * \nu) = \Phi(\mu) * \Phi(\nu)
\end{equation}
where the convolutions are in $\g$ and $G$, respectively.}
\end{theorem}

\noindent {\bf Remarks:}
\begin{list}{}{}
\item (i) Note that (\ref{wrapform}) implies the Duflo isomorphism for compact Lie groups since the Ad-invariant
distributions of support $\{0 \}$ in $\g$ are mapped by $\Phi$ to
central distributions of support $\{e \}$ in $G$.
\item (ii) A version of (\ref{wrapform}) has more recently been given by Andler, Sahi and Torossian (\cite{AST}) for all Lie groups which holds for {\it germs of
Ad-invariant hyperfunctions} with support at the identity (we leave
the precise details of the definition of germs of Ad-invariant
hyperfunctions to \cite{AST} and \cite{KV}). Their result was
conjectured by Kashiwara and Vergne (\cite{KV}), who proved it in
the case where $G$ is a solvable group.
\end{list}

This formula originated from the authors previous work on sums of
adjoint orbits (\cite{DRW}).  What (\ref{wrapform}) shows us is that
problems of convolution of central measures or distributions on a
(non-abelian) compact Lie group can be
transferred to Euclidean convolution of Ad-invariant distributions on $\g$.\\

Thus, since the solution to the Cauchy problem for the heat equation
can be written as a convolution between the heat kernel and the
initial data, we should be able to wrap the heat kernel on $\g \cong
\R^n$ to
that on $G$, and transfer the corresponding solution of the Cauchy problem.\\

Given the remarks in the Introduction, it is clearly of interest
also to consider whether there is a way to ``wrap Brownian motion'', and thus
obtain the heat kernel on $G$.  The main result of this paper will show how this is achieved \\

We will quickly recall the proof of Theorem \ref{thm2} here, since
it is instructive in our computation of the heat kernel on $G$.
Recall from section 2 the Kirillov character formula states that the
Fourier transform of the Liouville measure on the integral
co-adjoint orbit through $\lambda + \rho \in \Lambda^+$ is $j(X)
\chi_\lambda (\exp X)$.  The proof of Theorem \ref{thm2} follows
easily from this formula.  We give the elementary proof: Let $\pi
\in \widehat{G}$ have highest weight $\lambda \in \Lambda_+$, and
let $\mu^\wedge$ denote the Euclidean Fourier transform of $\mu$
with the convention
$$
\mu^\wedge (\xi) = \int_{\R^n} \mu (\x) e^{i \x \cdot \xi} d\x .
$$

Then we have:
\begin{lemma}\label{PhiHat}  Let $\mu$ be an Ad-invariant distribution of compact support on $\g$.  Then the Fourier transform of $\Phi(\mu)$ at $\pi$ is a multiple $c_\pi I_\pi$ of the identity, where
$$
c_\pi = (\Phi(\mu)) ^\wedge (\lambda + \rho) = \mu^\wedge (\lambda +
\rho) .
$$
\end{lemma}

{\bf Proof:}  By the definition of $\Phi$, $c_\pi$ is given by
$$
c_\pi = \frac{1}{d_\pi} \langle \Phi(\mu) , \chi_\lambda \rangle =
\frac{1}{d_\pi} \langle \mu , j \tilde{\chi}_\lambda \rangle .
$$

By applying Kirillov's character formula we have
$$
c_\pi = \langle \mu , \int_G e^{ig(\lambda + \rho)(\cdot)} dg
\rangle .
$$

By the $G$-invariance of $\mu$ this is
$$
c_\pi = \langle \mu , e^{i(\lambda + \rho)(\cdot)} \rangle =
\mu^\wedge (\lambda + \rho)
$$
as required. $\phantom{abcde} \square$\\

Theorem \ref{thm2} follows, since
\begin{align*}
(\Phi(\mu * \nu))^\wedge (\pi) & = (\mu * \nu)^\wedge (\lambda + \rho) I_\pi\\
& = \mu ^\wedge (\lambda + \rho) \nu ^\wedge (\lambda + \rho) I_\pi\\
& = (\Phi(\mu)) ^\wedge (\lambda + \rho) (\Phi(\nu)) ^\wedge (\lambda + \rho) I_\pi\\
& = (\Phi(\mu) * \Phi(\nu))^\wedge (\pi) . \phantom{abcde} \square
\end{align*}

\section{The wrap of Brownian motion}

We now give our main results on wrapping Brownian motion and heat
kernels on compact Lie groups.  Firstly, we show how to wrap the
Laplacian, the infinitesimal generator of a Brownian motion and the
heat semigroup.  We will then give the construction of how to ``wrap
Brownian motion".  Finally, by taking expectations we wrap the heat
kernel on $G$, thus providing an easy way to determine the
transition density of a Brownian motion on $G$ from that of $\g$.

\subsection{The wrap of the Laplacian}

From the Introduction, we note that one-half of the Laplacian
$$
\frac{1}{2} \Delta = \frac{1}{2} \sum_{i=1}^n
\frac{\partial^2}{\partial x_i^2}
$$
is the generator of associated heat semigroup $e^{t \Delta / 2}$ on
$\R^n$ and $G$.  In the next section, we will recall standard
results regarding Brownian motion on $\R^n$, where one-half of the
Laplacian is the infinitesimal generator of a Brownian motion on
$\R^n$.  Furthermore, one-half of the Laplacian on $G$ is the
infinitesimal generator of Brownian motion on $G$.\\

Thus, to see how to wrap a Brownian motion and the heat kernel from
$\g$ to $G$, we will first need to see how the infinitesimial
generator of the respective process and semigroup  - the Laplacian -
on $\g$ and $G$ are related by wrapping.  We will see that the
Laplacian on $\g$ is not quite mapped to the Laplacian on $G$.  More
precisely, we have:
\begin{prop}\label{wraplap}  Let $G$ be a compact connected Lie group with Lie algebra $\g$.  Let $L_G$ be the Laplacian on $G$ and $L_\g$ the Laplacian on $\g$.  Then for any $\mu \in
S_G (\g)$
$$
\Phi \bigl( L_\g (\mu) \bigr) = (L_G - \| \rho \|^2) \bigl( \Phi \mu
\bigr).
$$
\end{prop}

Firstly, we require the following:
\begin{lemma}\label{lemma332}  With the above notation, we have:\\
$$
\frac{\dim G}{24} = \| \rho \|^2 = -j^{-1} L_\g \, j
$$
\end{lemma}
\noindent {\bf Proof:}  $\frac{\dim G}{24} = \| \rho \|^2$ is
Freudenthal and de Vries' ``strange formula".  For the second
equality we use the Kirillov character formula.  Firstly, we need to
find the Fourier transform of $j$. Putting $\lambda = 0$ in the
Kirillov character formula we have that $j^\wedge (\xi) = \mu_\rho$,
where $\mu_\rho$ is the Liouville measure on $\mathcal{O}_\rho$.\\

Let $\nabla_\alpha$ be the directional derivative in the direction
of $\alpha$, so that the Laplacian is given by the gradient, that is
$\Delta = \nabla_\alpha \cdot \nabla_\alpha$.  By the elementary
properties of the Fourier transform, we have
$$
(\nabla_\alpha f)^\wedge (\lambda) = i \langle \alpha , \lambda
\rangle f^\wedge (\lambda).
$$

But $j^\wedge$ is supported on the co-adjoint orbit through $\rho$,
and so
$$
(\nabla_\alpha j)^\wedge (\rho) = i \langle \alpha , \rho \rangle
j^\wedge (\rho).
$$

Hence
$$ (\Delta j)^\wedge (\rho) = (\nabla_\alpha \cdot \nabla_\alpha j)^\wedge (\rho) = - \sum_{\alpha \in \Sigma^+} \langle \alpha , \rho \rangle \langle \alpha , \rho \rangle j^\wedge (\rho), $$
and thus
$$
(\Delta j)^\wedge (\rho) = - \| \rho \|^2 j^\wedge (\rho)
\phantom{abcde} \square
$$

\noindent {\bf Proof of Proposition \ref{wraplap}:} This is
essentially equivalent to \cite{HE3}, Ch. V, Proposition 5.1, where
it is proved for a more general class of symmetric spaces.  We will
give a proof of Proposition \ref{wraplap} for compact Lie groups
using the wrapping map: Let $\mu \in S_G (\g)$, Then,
\begin{align*}
\bigl\langle \Phi(L_\g \mu) & , f \bigr\rangle = \bigl\langle L_\g
\mu, j \tilde{f} \bigr\rangle \tag{by definition of the
wrap}\\
& = \bigl\langle \mu, L_\g (j \tilde{f}) \bigr\rangle \tag{since the Laplacian is a symmetric operator}\\
& = \bigl\langle \mu, j \, j^{-1} L_\g (j \tilde{f}) \bigr\rangle \\
& = \bigl\langle \mu, j . L_G^{\exp^{-1}} \tilde{f} + (L_\g j) \tilde{f} \bigr\rangle \tag{by \cite{HE2} Ch. II, Thm. 3.15}\\
& = \bigl\langle \mu, j . L_G^{\exp^{-1}} \tilde{f}\bigr\rangle +
\bigl\langle \mu, j \, j^{-1} (L_\g j) \tilde{f} \bigr
\rangle\\
& = \bigl\langle \mu, j . \widetilde{L_G f} \bigr\rangle +
\bigl\langle \mu, -j \| \rho \|^2 \tilde{f} \bigr\rangle \tag{ from
Lemma 2}\\
& = \bigl\langle \Phi(\mu), L_G f \bigr\rangle + \bigl\langle \Phi(\mu), -\| \rho \|^2 f \bigr\rangle\\
& = \bigl\langle L_G \Phi(\mu), f \bigr\rangle + \bigl\langle -\|
\rho \|^2 \Phi(\mu), f \bigr\rangle  \tag{since the Laplacian is a symmetric operator}\\
& = \bigl\langle (L_G - \| \rho \|^2)(\Phi(\mu)), f \bigr\rangle
\phantom{abcde} \square
\end{align*}

{\bf Remark:} $L_G - \| \rho \|^2$ is also known as the {\bf shifted
Laplacian}.  It has often been observed that it is more appropriate
in aspects of harmonic analysis.  As an example, it is shown in
\cite{HE3}, Ch. V, $\S 5$ that Huyghen's principle does not hold for
a compact Lie group, but does (in odd dimension) when the shifted
Laplacian is used.  We will also refer to Brownian motion with a
potential of $\| \rho \|^2$ as a {\bf shifted Brownian
motion} in accordance with its generator, the shifted Laplacian.\\

\subsection{The wrap of the It\^o equation}

Proposition \ref{wraplap} shows that $\Phi \bigl( L_\g (u) \bigr) =
(L_G - \| \rho \|^2) \bigl( \Phi u \bigr)$.  This formula allows us
to ``wrap'' the Laplacian from $\g$, to the
shifted one on $G$.\\

The actual mechanics of {\it wrapping Brownian motion} are not
immediately obvious.  The wrapping map is a homomorphism from the
algebra of Ad-invariant distributions on $C^\infty (\g)$ to the
algebra of central distributions on $C^\infty (G)$, defined by $\phi
\mapsto \phi \iota$ where $\iota : f \mapsto j.f \circ \exp$. In
this subsection we will show to wrap Brownian motion by using stochastic differential equations.\\

We wish to emphasise at this point that wrapping Brownian motion is
not to be interpreted as a map $\Phi : \zeta_t \rightarrow \xi_t$ of
Brownian motion on $\R^n$ to $G$ (such as the well-known It\^{o}
map), but rather it is a method of ``converting'' between the S.D.E.'s of these two processes.\\

%%%%%%%%%%%%%%%%%%%%%%%%%%%%%%%%%%%%%%%%%%%%%%%%%%%%%%%%%%%%%%%%
%
%  BM on R^n
%
%%%%%%%%%%%%%%%%%%%%%%%%%%%%%%%%%%%%%%%%%%%%%%%%%%%%%%%%%%%%%%%%

We recall the definitions regarding Brownian motion and stochastic
integration on $\R^n$ from \cite{IKW}, \cite{KUN} and \cite{RY}.  We
use the standard notations regarding probability spaces, filtrations
and expectations from these sources.  We briefly state the following
definitions for the purposes of notation.

\begin{df} A {\bf (standard) Brownian motion} on $\R$ is a continuous stochastic process $(B_t)_{t \geq 0}$ such that for $0 \leq s < t < \infty$:
\begin{list}{}{}
\item (i) For $ 0 \leq s < t < \infty$, $B_t - B_s$ is a normally distributed random variable with mean 0 and variance
$t-s$.
\item (ii) For $ 0 \leq t_0 < t_1 < \dots < t_n < \infty$,
$\{ B_{t_k} - B_{t_{k-1}}, \, k = 1 , \dots , n \}$ is a set of
independent random variables.
\end{list}

Furthermore, an {\bf $n$-dimensional (standard) Brownian motion} on
$\R^n$ is a right continuous stochastic process
$$
B_t = (B_t^{(1)}, \dots, B_t^{(n)})
$$
where each $\{B_t^{(i)}\}_{t \geq 0} $ is an independent Brownian
motion.

\end{df}

It\^o's theory of stochastic integration provides us with the
following formula:
\begin{theorem}\label{MDIF} (Multidimensional It\^o formula) (\cite{CHW} Thm. 5.10)  Let $(M_t)_{t \geq 0}$ be a continuous local martingale with values in $\R^n$.  Suppose $f$ is a continuous $C^2$ function $f: \R^n \times \R^+ \rightarrow \R$.  Then a.s. for each $t > 0$,
\begin{align*}
f(M_t, t) - f(M_t, 0) = & \sum_{i=1}^n \int_0^t \frac{\partial f}{\partial x_i} (M_s,s) dM^{(i)}_s + \int_0^t \frac{\partial f}{\partial t} (M_s, s) ds\\
& + \sum_{i=1}^n \sum_{j=1}^n \int_0^t \frac{\partial^2 f}{\partial
x_i \partial x_j} (M_s, s) d \langle M_s^{(i)}, M_s^{(j)} \rangle_s
\end{align*}
\end{theorem}

Stratonovich developed a stochastic integral that would ``conform''
to the usual rules of calculus, known as the {\it Stratonovich
integral}:
\begin{df}  Let $B$ be a Brownian motion and $X$ an $L^2$-martingale.  The {\bf Stratonovich integral}, denoted by $\int_0^t X_s \circ dB_s$, is defined by
\begin{equation}\label{strat}
\int_0^t X_s \circ dB_s = \int_0^t X_s dB_s + \tfrac 12 \langle X, B
\rangle_t
\end{equation}
\end{df}

From these formula, we may construct Brownian motion, $(\zeta_t)_{t
\geq 0}$, on $\R^n$ as the solution to the Stratonovich S.D.E.:
\begin{equation}
h(\zeta_t) = h(0) + \sum_{i=1}^n \int_0^t \frac{\partial h}{\partial
x_i} (\zeta_s) \circ dB_s^{(i)} , \phantom{abcde} \zeta_0 = 0
\end{equation}
or, equivalently as the It\^o S.D.E.:
\begin{equation}
h(\zeta_t) = h(0) + \sum_{i=1}^n \int_0^t \frac{\partial h}{\partial
x_i} (\zeta_s) dB_s^{(i)} + \tfrac 12 \sum_{i=1}^n \int_0^t
\frac{\partial^2 h}{\partial x_i^2} (\zeta_s) ds, \phantom{abcde}
\zeta_0 = 0
\end{equation}
for any $h \in C^\infty_b(\R^n)$.  Thus, the Laplacian
$\frac{\partial^2 h}{\partial x_i^2}$, is said to be the
infinitesimal generator of Brownian motion.\\

%%%%%%%%%%%%%%%%%%%%%%%%%%%%%%%%%%%%%%%%%%%%%%%%%%%%%%%%%%%%%%%%
%
%  BM on Lie Groups
%
%%%%%%%%%%%%%%%%%%%%%%%%%%%%%%%%%%%%%%%%%%%%%%%%%%%%%%%%%%%%%%%%

We now define Brownian motion on a semisimple Lie group:
\begin{df}  Suppose $G$ is an $n$-dimensional semisimple Lie group with Lie algebra $\g$, and $X_1 ,\dots, X_m$ are vector fields on $G$.  If $(B_t)_{t \geq 0}$ is an $m$-dimensional Brownian motion on $\g$ and $p \in G$, then a {\bf $G$-valued stochastic process} $(\xi_t)_{t \geq 0}$ is said to be a {\bf solution of}
\begin{equation}
d\xi_t = \sum_{i=1}^n X_i (\xi_t) \circ dB_t^{(i)}, \phantom{abcde}
\xi_0 = p
\end{equation}
if for each $f \in C^\infty (G)$ we have
\begin{equation}\label{232}
f(\xi_t) = f(p) + \sum_{i=1}^n \int_0^t (X_i f)(\xi_s) \circ
dB_s^{(i)}
\end{equation}

The solution of (\ref{232}) is a {\bf Brownian motion} on $G$,
starting at $p$.
\end{df}

The above definition is often given (\cite{TAY}, \cite{LIAO2}) for
Brownian motion on semisimple Lie groups.  How this equation arises
from the construction of Brownian motion by ``rolling without
slipping" on Riemannian structure (\cite{IKW}) on a semisimple Lie
group is given in \cite{COR}, pp. 64-66. Converting Stratonovich to
It\^o integrals in the corresponding integral equation (\ref{232})
(see \cite{LIAO2}, Ch. 1) yields:
\begin{equation}\label{233}
f(\xi_t) = f(e) + \sum_{i=1}^n \int_0^t (X_i f)(\xi_s) dB_s^{(i)} +
\frac{1}{2} \sum_{i=1}^n \int_0^t (X_i^2 f)(\xi_s) ds
\end{equation}

(\ref{233}) implies that the generator of Brownian motion on $G$ is
given by one half of the Laplacian. This is the definition of a
Brownian motion given in \cite{SHG}.  We also refer the reader to
\cite{PAU} for a comparison of these two constructions.\\

Following from the Introduction, we have that
$$
(P_t f)(g) = \E_g (f(\xi_t)) = \E(f(g \cdot \xi_t))
$$
and
$$
(P_t f)(g) = \E \Bigl( \int_0^t (\tfrac 12 L_G f)(g \cdot \xi_s)ds
\Bigr)
$$

That is, the heat kernel is the law of Brownian motion, and is also
the fundamental solution of the heat semigroup.  It also follows
that the infinitesimal generator of the heat semigroup is equal to
the Laplacian on $C_c (G)$.  Furthermore,
$$
P(g \cdot \xi_t \in dh ) = p_t (g^{-1} h) dh
$$
and therefore
$$
\E_g(f(\xi_t)) = \E(f(g \cdot \xi_t) = \int_G f(h) p_t(g^{-1} h) dh
$$

%%%%%%%%%%%%%%%%%%%%%%%%%%%%%%%%%%%%%%%%%%%%%%%%%%%%%%%%%%%%%%%%
%
%  Wrapping BM
%
%%%%%%%%%%%%%%%%%%%%%%%%%%%%%%%%%%%%%%%%%%%%%%%%%%%%%%%%%%%%%%%%

We now consider the mechanics of wrapping Brownian motion. View $\g$
as $\R^n$ and consider Brownian motion, $(B_t)_{t \geq 0}$, on $\g$.
Hence, we may regard $(B_t)_{t \geq 0}$ as a process on the tangent
spaces of both $\g$ and $G$. As before, we may construct Brownian
motion $(\zeta_t)_{t \geq 0}$ on $\g$ as the solution to the S.D.E.:
\begin{equation}\label{Itoong}
h(\zeta_t) = h(0) + \sum_{i=1}^n \int_0^t \frac{\partial h}{\partial
x_i} (\zeta_s) dB_s^{(i)} + \tfrac 12 \sum_{i=1}^n \int_0^t
\frac{\partial^2 h}{\partial x_i^2} (\zeta_s) ds, \phantom{abcde}
\zeta_0 = 0
\end{equation}
for any $h \in C^\infty_b(\g)$.\\

This is a kind of ``identity S.D.E.", where the solution is $B_t$.
Likewise, we could na\"ively define ``shifted Brownian motion"  -
corresponding to the shifted Laplacian  -  on $G$ as the solution to
the S.D.E.:
\begin{equation}\label{ItoonG}
f(\xi_t) = f(e) + \sum_{i=1}^n \int_0^t (X_i f)(\xi_s) dB_s^{(i)} +
\tfrac 12 \sum_{i=1}^n \int_0^t ((X_i^2 - \|\rho\|^2) f)(\xi_s)ds,
\phantom{abcde} \xi_0 = e.
\end{equation}
where $\bigl( X_i \bigr)_{i=1}^n$ is an orthonormal basis of the Lie
algebra.\\

However, there is a problem with this definition: $\xi_t$ is
generated by $L_G - \| \rho \|^2$ and is thus not Markovian, which
is contrary to it being a solution of \ref{ItoonG}.  We thus define
$\xi_t$ by starting with a standard Brownian motion, $\tilde{\xi}_t$
on $G$, and ``killing" the process by applying the Feynman-Ka\v{c}
Theorem to obtain shifted Brownian motion on the Lie group.  Note
that this also involves enlarging the probability space from
$\Omega$ to $\Omega \times [0,T]$, equipped
with the appropriate product measure, $P$.\\

By writing shifted Brownian motion as a solution to an SDE we have:
\begin{lemma}\label{FKlemma}
Suppose $(\tilde{\xi}_t)_{t \geq 0}$ is the solution on the filtered
probability space $(\Omega \times [0,T], \mathcal{F}, \mathcal{F}_t,
P)$ to the SDE
$$
f(\tilde{\xi}_t) = f(e) + \sum_{i=1}^n \int_0^t (X_i
f)(\tilde{\xi}_s) dB_s^{(i)} + \tfrac 12 \sum_{i=1}^n \int_0^t
(X_i^2 f)(\tilde{\xi}_s)ds
$$
Consider a new measure of the form
$$
d\tilde{P} = e^{-C t} dP,
$$
where $C > 0$ is a constant, and $t \in [0,T]$.  Then $(\xi_t)_{t
\geq 0}$ is the solution on the filtered probability space $(\Omega
\times [0,T], \mathcal{F}, \mathcal{F}_t, \tilde{P})$ to the SDE
$$
f(\xi_t) = f(e) + \sum_{i=1}^n \int_0^t (X_i f)(\xi_s) dB_s^{(i)} +
\tfrac 12 \sum_{i=1}^n \int_0^t ((X_i^2 - C)f)(\xi_s)ds.
$$
\end{lemma}

\begin{df}\label{shiftedBM}  We refer to the solution of
$$
f(\xi_t) = f(e) + \sum_{i=1}^n \int_0^t (X_i f)(\xi_s) dB_s^{(i)} +
\tfrac 12 \sum_{i=1}^n \int_0^t ((X_i^2 - C)f)(\xi_s)ds.
$$
on the filtered probability space $(\Omega \times [0,T],
\mathcal{F}, \mathcal{F}_t, \tilde{P})$ as a {\bf shifted Brownian
motion}.
\end{df}

We now define what it means to {\it wrap Brownian motion}
\begin{df}\label{wrapbm}  Let $\zeta_t$ be a Brownian motion on $\g \cong \R^n$.  The wrap of $\zeta_t$ is to a process $\xi_t$ on $G$ is given by the mapping $f \mapsto j.f \circ \exp = h$, where $\xi_t$ is the shifted Brownian motion in Definition \ref{shiftedBM}.  We write this as
$$
\Phi (\zeta_t) = \xi_t
$$
\end{df}

Furthermore, it can be seen that $\xi_t$ is a Brownian motion on $G$
with a potential of $\| \rho \|^2$:  By replacing $f \in C^\infty
(G)$ with $j.f \circ \exp \in C^\infty_c (\g)$, (\ref{ItoonG})
becomes
\begin{align*}
j.f(\exp(\xi_t)) & = j.f(\exp(e)) + \sum_{i=1}^n \int_0^t (X_i j.f)(\exp(\xi_s)) dB_s^{(i)}\\
&\phantom{abcde} + \tfrac 12 \sum_{i=1}^n \int_0^t (X_i^2
j.f)(\exp(\xi_s))ds - \tfrac 12 \|\rho\|^2 \int_0^t j.f(\exp(\xi_s))
ds.
\end{align*}

Letting $j.f \circ \exp = h \in C^\infty_0 (\g)$, we obtain
\begin{align*}
h(\zeta_t) & = h(0) + \sum_{i=1}^n \int_0^t (X_i (j \tilde f))(\zeta_s) dB_s^{(i)} + \tfrac 12 \sum_{i=1}^n \int_0^t (X_i^2 j \tilde f)(\zeta_s)ds - \tfrac 12 \int_0^t \|\rho\|^2 j \tilde f(\zeta_s) ds\\
& = h(0) + \sum_{i=1}^n \int_0^t \frac{\partial h}{\partial x_i}(\zeta_s) dB_s^{(i)} + \tfrac 12 \int_0^t (L_G^{\exp^{-1}} h)(\zeta_s)ds - \tfrac 12 \int_0^t \|\rho\|^2 h(\zeta_s) ds\\
& = h(0) + \sum_{i=1}^n \int_0^t \frac{\partial h}{\partial x_i}(\zeta_s) dB_s^{(i)} + \tfrac 12 \int_0^t (L_\g h)(\zeta_s)ds\\
&\phantom{aaaaabbbbbcccccdddddeeeee} + \tfrac 12 \int_0^t \|\rho\|^2 h(\zeta_s) ds - \tfrac 12 \int_0^t \|\rho\|^2 h(\zeta_s) ds\\
& = h(0) + \sum_{i=1}^n \int_0^t \frac{\partial h}{\partial
x_i}(\zeta_s) dB_s^{(i)} + \tfrac 12 \sum_{i=1}^n \int_0^t
\frac{\partial^2 h}{\partial x_i^2}(\zeta_s)ds,
\end{align*}
which is (\ref{Itoong}).\\

We now show that taking expectations of our Brownian motion on $\g$
and wrapped Brownian motion on $G$ corresponds to
wrapping the heat kernel from $\g$ to $G$.\\

%Proposition
Following from Definition \ref{wrapbm}, the wrap of Brownian motion
on $\g$ is a shifted Brownian motion on $G$ which satisfied the
S.D.E.
\begin{equation}\label{345}
f(\xi_t) = f(e) + \sum_{i=1}^n \int_0^t (X_i f)(\xi_s) dB_s^{(i)} +
\tfrac 12 \sum_{i=1}^n \int_0^t (X_i^2 f)(\xi_s)ds - \tfrac 12
\|\rho\|^2 \int_0^t f(\xi_s) ds
\end{equation}

Since the $X_i$'s are left-invariant, the process starting at $g \in
G$ is $(g \cdot \xi_t)_{t \geq 0}$.  We consider the law of this
process, and define the operators $P_t$ on $C_c (G) \subseteq C_0
(G)$ by
$$
(P_t f)(g) = \E_g^{\tilde{P}} (f(\xi_t)) = \E^{\tilde{P}}(f(g \cdot
\xi_t)).
$$

By (\ref{345}) we have
$$
f(g \cdot \xi_t) = f(g) + \sum_{i=1}^n \int_0^t (X_i f)(g \cdot
\xi_s) dB_s^{(i)} + \tfrac 12 \sum_{i=1}^n \int_0^t (X_i^2 f)(g
\cdot \xi_s)ds - \tfrac 12 \|\rho\|^2 \int_0^t f(g \cdot \xi_s) ds.
$$
under the measure $\tilde{P}$. It now follows that
\begin{align*}
(P_t f)(g) & = \E_g^{\tilde{P}} (f(\xi_t)) = \E^{\tilde{P}}(f(g \cdot \xi_t))\\
& = \E^{\tilde{P}} \Bigl( \tfrac 12 \sum_{i=1}^n \int_0^t (X_i^2
f)(g \cdot \xi_s)ds - \tfrac 12 \|\rho\|^2 \int_0^t f(g \cdot \xi_s)
ds \Bigr)\\
& = e^{-Ct} \E^P \Bigl( \tfrac 12 \sum_{i=1}^n \int_0^t (X_i^2 f)(g
\cdot \tilde{\xi}_s)ds \Bigr)
\end{align*}

Furthermore, if for each $t$ the random variable $\xi_t$ has a
density $p_t$ with respect to the Haar measure $dh$ then
$$
\E^P_g(f(\xi_t)) = \E^P (f(g \cdot \tilde{\xi}_t)) = \int_G f(h)
p_t(g^{-1} \cdot h) dh.
$$

To calculate the wrap, we need to replace $f$ by $j. f \circ \exp =
h$.  On $\g$ this gives
\begin{align*}
(Q_t (j. f \circ \exp))(X) = (Q_t h)(X) & = \E_X(h(\zeta_t)) = \E(h(X + \zeta_t))\\
& = \E \Bigl( \tfrac 12 \sum_{i=1}^n \int_0^t \frac{\partial^2
h}{\partial x_i^2} (X + \zeta_t) dt \Bigr).
\end{align*}

Thus we may write
$$
\E_X(j. f \circ \exp(\zeta_t)) = \E^{\tilde{P}}_{\exp X} (f(\xi_t)).
$$

In terms of the semigroup of operators for the heat equation we have
$$
P_t f = Q_t (j. f \circ \exp).
$$
Using the wrapping formula, we may rewrite this as
$$
\Phi(P_t) = Q_t.
$$

We have proved:
\begin{theorem}\label{wrapBM2}  Suppose $\xi_t$ is the wrap of the Brownian motion $\zeta_t$ on $\g$.  Then the law of $\xi_t$ may be found by wrapping the law of Brownian motion on its Lie algebra.  That is,
$$
\E_X(j \cdot f \circ \exp(\zeta_t)) = \E^{\tilde{P}}_{\exp X}
(f(\xi_t)).
$$

The transition density of $\xi_t$ is given by
$$
\Phi (p_t)(\exp H) = q_t^\rho (g)
$$
where $p_t(x)$ is the heat kernel on $\g = \R^n$, and $q_t^\rho (g)$
is the heat kernel corresponding to the shifted Laplacian on $G$
\end{theorem}

Let $(\xi_t)_{t \geq 0}$ be wrapped Brownian motion from
(\ref{345}).  Then the expectation of $(\xi_t)_{t \geq 0}$ is the
shifted heat kernel:
\begin{equation}
\E^{\tilde{P}}(\xi_t) = q^\rho_t (g) = \sum_{\lambda \in \Lambda^+}
d_\lambda \chi_\lambda (g) e^{-\| \lambda + \rho \|^2 t/2} ,
\phantom{abc} t \in \R^+, \; g \in G.
\end{equation}

This expectation is taken with respect to the measure $P$.  By
applying Lemma \ref{FKlemma} (with $C = \tfrac 12 \| \rho \|^2$) to
Theorem \ref{wrapBM2} we get
\begin{prop}
The expectation of $(\xi_t)_{t \geq 0}$ under $P$, is
$$
\E^P(\xi_t) = \sum_{\lambda \in \Lambda^+} d_\lambda \chi_\lambda
(g) e^{-(\| \lambda + \rho \|^2 - \|\rho \|^2)t/2} , \phantom{abc} t
\in \R^+, \; g \in G.
$$
\end{prop}

\noindent {\bf Proof:}  Taking expectations under $\tilde{P}$ yields
\begin{align*}
\E^P (\xi_t) & = \E^P (e^{\| \rho \|^2 t/2} \xi_t) = e^{\| \rho \|^2 t/2} \E^{\tilde{P}} (\xi_t)\\
& = e^{\| \rho \|^2 t/2}q_t^\rho (g) = q_t (g)\\
& = \sum_{\lambda \in \Lambda^+} d_\lambda \chi_\lambda (g) e^{-(\|
\lambda + \rho \|^2 - \|\rho \|^2)t/2} , \phantom{abc} t \in \R^+,
\; g \in G.
\end{align*}
as required. $\square$\\

This generalises the situation of $\R^n$ as shown in the
Introduction.

\subsection{The wrap of the heat kernel}

Let $p_t(x)$ be the heat kernel on $\R^n$ be given by
\begin{equation}
p_t(x) = (2\pi t)^{-n/2} e^{-\frac{\| \x \|^2}{2t}}, \phantom{abcde}
t \in \R^+, \; \x \in \R^n.
\end{equation}
and $q_t (g)$ is the heat kernel on $G$ be given by

\begin{equation}
q_t (g) = \sum_{\lambda \in \Lambda^+ } d_\lambda \chi_\lambda (g)
e^{-c(\lambda) t/2} , \phantom{abcde} t \in \R^+, \; g \in G.
\end{equation}
where $c(\lambda) = \| \lambda + \rho \|^2 - \| \rho \|^2$ (Prop.
2.7.20). We write the heat kernel corresponding to the shifted
Laplacian on $G$ as $q_t^\rho (g)$. It is given by
\begin{equation}
q_t^\rho (g) = \sum_{\lambda \in \Lambda^+} d_\lambda \chi_\lambda
(g) e^{-\| \lambda + \rho \|^2 t/2} , \phantom{abc} t \in \R^+, \; g
\in G.
\end{equation}

We now use Theorem \ref{wrapBM2} to derive the heat kernel on $G$ by
calculating $\Phi (p_t)$.  In fact, we obtain a slightly stronger
version of this formula for $G$-invariant Schwartz functions:
\begin{prop}\label{wrapS} Let $\mu \in S(\g)$ be $G$-invariant, and $\hat{\mu}$ its (Euclidean) Fourier transform.  Then $\Phi(\mu) \in C^\infty_G (G)$ is given on $T$ by
$$
\Phi(\mu) (t) = \sum_{\lambda \in \Lambda^+} d_\lambda \, \hat{\mu}
(\lambda + \rho) \chi_\lambda (t), \phantom{abcde} \forall t \in T.
$$
\end{prop}
\noindent {\bf Proof:}  From Lemma \ref{PhiHat} we have that
$$
\Phi^\wedge (\mu) (\pi_\lambda)  = \hat{\mu} (\lambda + \rho)
I_{\pi_\lambda}
$$
and we invert the Fourier transform to obtain
\begin{equation}
\Phi(\mu) (t) = \sum_{\lambda \in \Lambda^+} d_\lambda \, \hat{\mu}
(\lambda + \rho) \chi_\lambda (t), \phantom{abcde} \forall t \in T.
\end{equation}
as required. $\phantom{ab} \square$\\

The shifted heat kernel on $G$ may be calculated by wrapping the
heat kernel on $\g$ using Proposition \ref{wrapS} and Theorem
\ref{thm1}. This in turn recovers the formulae of Sugiura
(\cite{SUG}) and Eskin (\cite{ESK}) for the heat kernel on a compact
Lie group. We also note that our expressions for the heat kernel can
also been seen in \cite{URA} and \cite{CMR}, which use the Poisson
summation formula.
\begin{theorem}\label{wrapHK} Let $p_t(x) = (2\pi t)^{-n/2} e^{-\frac{\| \x \|^2}{2t}}, \; t \in \R^+, \; \x \in \g$ be the heat kernel on $\g$.  Then $\Phi(p_t)$ is the shifted heat kernel on $G$ (given on $T$), given by
\begin{align}
\Phi(p_t) (\exp H) & = \sum_{\lambda \in \Lambda^+} d_\lambda \, e^{-\| \lambda + \rho \|^2 t} \chi_\lambda (\exp H)\label{356}\\
& = (2\pi t)^{-n/2} \sum_{\gamma \in \Gamma} e^{\frac{-\|H + \gamma
\|^2}{2t}} \frac{1}{j(H + \gamma)}\label{357}
\end{align}
for all $H \in \tg.$
\end{theorem}

\noindent {\bf Proof of Theorem \ref{wrapHK}:}   Setting $\mu = p_t$
in Proposition \ref{wrapS} gives us our result for wrapping the heat
kernel:
$$
\hat{p}_t(\xi) = e^{-\| \xi \|^2 t/2}
$$
and therefore
$$
\Phi^\wedge (p_t) (\pi_\lambda) = e^{- \| \lambda + \rho \|^2 t/2}.
$$

Thus by Proposition \ref{wrapS} we have
$$
\Phi(p_t) (\exp H) = \sum_{\lambda \in \Lambda^+} d_\lambda \,
e^{-\| \lambda + \rho \|^2 t} \chi_\lambda (\exp H), \phantom{abcde}
\forall H \in \tg.
$$
giving (\ref{356}), which is the heat kernel corresponding to the shifted Laplacian.\\

By Theorem \ref{thm1} we may wrap the heat kernel $p_t$ by putting
$$
\Phi (p_t)(\exp H) = \Phi (j (p_t \tfrac 1j ))(\exp H) = (2\pi
t)^{-n/2} \sum_{\gamma \in \Gamma} e^{\frac{-\|H + \gamma \|^2}{2t}}
\frac{1}{j(H + \gamma)}
$$
which yields (\ref{357}).  This is valid for the regular points of
$G$. It is also valid for the singular points since (\ref{356}) is
$C^\infty$, and therefore (\ref{357}) is also $C^\infty$ since it is
clearly $C^\infty$ on the regular elements, and is thus $C^\infty$
at the singular points by analytic continuation.\\

Proposition \ref{wrapS} also allows us to prove Proposition
\ref{wraplap} by considering the Laplacian as a distribution
supported at the identity, acting by convolution. Write the
Laplacian on $\g$ as a Fourier multiplier:
$$
\widehat {L_\g f} (\xi) = -\| \xi \|^2 \hat{f} (\xi)
$$

Now, $ \Phi^\wedge (\mu)  = \hat{\mu} (\lambda + \rho)$, so by
taking $\hat{\mu}(\xi) = -\|\xi\|^2$,
\begin{equation}\label{lapmult}
\Phi(L_\g) (t) = -\sum_{\lambda \in \Lambda^+} d_\lambda \, \|
\lambda + \rho \|^2 \chi_\lambda (t), \phantom{abcde} \forall t \in
T.
\end{equation}
which is the shifted Laplacian on a compact Lie group, given by a
distribution supported at the identity. However, the Laplacian as a
distribution supported at the identity is given by
\begin{equation}\label{lapmult2}
(L_G) (t) = -\sum_{\lambda \in \Lambda^+} (\| \lambda + \rho \|^2 -
\| \rho \|^2) \chi_\lambda (t)
\end{equation}

The discrepancy between \ref{lapmult} and \ref{lapmult2} yields the
following (c.f. Proposition \ref{wraplap}):
$$
\Phi \bigl( L_\g \bigr) = (L_G - \| \rho \|^2) \bigl( \Phi \bigr)
$$

\noindent {\bf Remarks:} The heat kernel on a compact Lie group has
been previously derived by Arede \cite{ARE} and Watanabe \cite{WAT},
although their methods are quite different from ours.  We briefly
summarise their results here:\\

In \cite{ARE}, the formula for the heat kernel on a compact Lie
group is given by
\begin{equation}\label{361}
q_t (\exp H) =  (2\pi t)^{-d/2} j(H)^{-1} e^{\frac{\|H\|^2}{2t} +
\|\rho\|^2 t/2} E(\chi_{\tau > t})
\end{equation}
where $\chi_{\tau > t}$ is the indicator function of the first exit
time of the so-called ``Brownian Bridge'' from the fundamental
domain.  Arede's proof involves the Elworthy-Truman ``Elementary
Formula''.  The heat kernel for the group $SU(2)$ is then given as
\begin{equation}
q_t (g) =  (2\pi t)^{-3/2} \sum_{j \in \Z} \frac{4\sqrt{2} j \pi +
|\lambda|} {2\sqrt{2} \sin \bigl[ (4\sqrt{2} j \pi + |\lambda|)
(2\sqrt{2}) \bigr]} e^{|4\sqrt{2}j\pi + \lambda|^2/2t} e^{t/16}
\end{equation}
where $\lambda \in \R$ is such that $|\lambda| = d(g,e)$ and $|\lambda| < 2\sqrt{2} \pi$.\\

This is similar to the formula given in \cite{WAT} for the group
$SU(2)$, but with different normalisations.  Watanabe's formula is:
\begin{equation}
q_t (\exp H) =  (2\pi t)^{-3/2} \exp \biggl\{ \frac{1}{4} t \biggr\}
\sum_{n \in \Z} e^{\frac{(H + n)^2}{2t}} \frac{ \frac{H + n}{2} }{
\sin \bigl( \frac{H + n}{2}\bigr)}, \;\;\; H \in \R
\end{equation}

Watanabe exploits the fact that the Laplace transform of the L\'evy
stochastic area process is $j^{-1}$.  Both Arede's and Watanabe's
work can be derived from a general formula on Riemannian manifolds
known as the {\it Minakushisudarum-Pleijel expansion}, which we will
examine in the sequel to this paper.

\subsection{Remarks on other processes}

In this section we show how the wrapping map can be used to transfer
other stochastic processes from $\g$ to $G$.  We firstly consider
the results of Kingman on spherically symmetric random walks
(\cite{KIN}) and show how in the case of $\R^3$ how the wrapping map
naturally transfers these to random walks on the conjugacy classes
of $SU(2)$. We then use the wrapping map to deduce certain recent
results of Liao (\cite{LIAO}) on the distribution of
$G$-invariant {\it L\'evy process}.\\

%, and the Central Limit Theorem of Clerc and Roynette (\cite{CR}).\\

%\subsection{Symmetric random walks}

We now consider spherically symmetric random walks studied by
Kingman \cite{KIN}: take two independent random variables $X$ and
$Y$ in $\R^3$, with lengths $|X|$ and $|Y|$, but with their
direction uniformly distributed.  The sum $Z = X+Y$ is uniformly
distributed in terms of its direction, but its length
$|Z|$ is a random number in the range $|X-Y| \leq |Z| \leq |X+Y|$.\\

In general, if the probability distributions of $|X|$ and $|Y|$ are
$\mu_X, \, \mu_Y \in M_1 (\R^+)$ (respectively), then $|Z|$ is a
random variable with probability distribution $\mu_Z$ depending on
$\mu_X$ and $\mu_Y$, with $\mu_Z = \mu_X * \mu_Y$.  This is
precisely the situation relating to the sums of adjoint orbits
considered in \cite{DRW}.\\

{\bf Remark:} The structures of Adjoint obits on $\g$, and conjugacy
classes in $G$ form structures known as {\it hypergroups} under the
operation of convolution  -  the reader is referred to \cite{WLD}
for further details.  The wrapping map forms an algebra isomorphism
between these two structures.\\

We now consider the wrapping map.  From the wrapping formula we have
that
\begin{equation}\label{341}
\Phi(\mu *_\g \nu) = \Phi(\mu) *_G \Phi(\nu)
\end{equation}

Recall that the Adjoint orbits in $\g$ are mapped to conjugacy
classes in $G$ by the relation $C_i = \exp \mathcal{O}_i$ via the
formula $\exp \Ad(g) X = g^{-1} \exp X g$. As a consequence of
(\ref{341}) we have:
\begin{prop}  Suppose $X$ and $Y$ are spherically symmetric random walks in $\R^3$, with the probability distributions of $|X|$ and $|Y|$ being $\mu_X$ and $\mu_Y$ (respectively), then the distribution of the wrap of $|X +
Y|$ on $SU(2)$ is
\begin{equation}
\Phi(\mu_X * \mu_Y) = \Phi(\mu_X) * \Phi(\mu_Y)
\end{equation}
\end{prop}

Also recall from Lemma \ref{PhiHat} that $\Phi(\mu)^\wedge
(\pi_\lambda) = \mu^\wedge (\lambda + \rho) I_\lambda$.\\

Following the introduction in \cite{KIN}, the characteristic
function of a spherically symmetric random walk on $\R^n$ is given
by
$$
\phi_X (\mathbf{t}) = \E(e^{i \mathbf{t} X}) = E(e^{i t X \cos
\theta})
$$
where $t = \| \mathbf{t} \|$, and $\theta$ is the angle between the
vectors $\mathbf{t}$ and $X$.  It is then shown by Kingman that
$$
E(e^{i x \cos \theta})  = \frac{J_{(n/2) - 1} (x)}{(x/2)^{(n/2) -
1}} ((n/2) - 1)!
$$
where $J_\lambda (\cdot)$ is the Bessel function of the first kind
of order $\lambda$, given by
\begin{equation}\label{3323}
\frac{J_{(n/2) - 1} (\lambda x)}{(\lambda x/2)^{(n/2) - 1}} =
\int_{S^{n-1}} e^{i\lambda \langle x,\omega \rangle} d\omega.
\end{equation}

Here, the Riemannian measure $d\omega$ has mass $((n/2) - 1)!$. Note
that \ref{3323} is the Kirillov character formula for a
compact Lie group, given in terms of generalised Bessel functions.\\

Let $n = 2(1 + s)$ and let $\Lambda_s (x) = J_s (x) s! (x/2)^{-s}$.
We have
$$
\phi_X (\mathbf{t}) = \E(\Lambda_s (tX))
$$

Since we are considering independent, spherically symmetric random
vectors, we will use the radial characteristic function
$$
\Psi_X (\mathbf{t}) = \E(\Lambda_s (tX))
$$
which in the case of $\R^3$ is
\begin{align*}
\Psi_X (\mathbf{t}) & = \int_0^\infty  \mu_X (x) \Lambda_s (tx) dx \\
& = \int_0^\infty  \mu_X (x) \int_{S^{2}} e^{it(x,\omega)} d\omega \, dx \\
& = \int_0^\infty  \mu_X (x) e^{it(x)} dx
\end{align*}
so that we have
\begin{equation}\label{343}
\Psi_X (\lambda + \rho) = \mu_X^\wedge (\lambda + \rho) =
\Phi^\wedge (\mu_X) (\pi)
\end{equation}

(\ref{343}) may then be inverted to obtain the transition density of
the random walk on $SU(2)$.  We now generalise some of our results
on Brownian motion to other processes using the above results.  In
this section we will consider {\it L\'evy processes}:
\begin{df} A {\bf L\'evy process}, $g_t$, is a stochastic process with independent and stationary increments, which has right continuous paths with left hand limits.
\end{df}

This definition includes both discrete and continuous processes. For
further details on L\'evy processes on Lie groups, the reader is
referred to \cite{LIAO} and \cite{LIAO2}.  We will assume these processes to start at the identity $e$ in $G$.\\

In analogy with the case of where the Laplacian is the generator of
Brownian motion and heat transition semigroup, it is also well known
that L\'evy processes have a {\it Feller transition semigroup},
$e^{t \mathcal{L} /2}$, with generator $\mathcal{L}$ that gives rise
to a unique semigroup of convolution operators $P_t$ which may be
convolved with the initial data $f(x) = u(0,x)$ to give the
transition density:
\begin{align*}
u(x,t) & = e^{t \mathcal{L} / 2} f(x)\\
& = P_t f(x)\\
& = (p_t * f)(x)\\
& = \int_G p_t(x^{-1} y) f(y) dy .
\end{align*}

Similarly, for any $f \in C^\infty (G)$ the distribution of $g_t$ is
completely determined by its generator, $\mathcal{L}$, given by
$$
\mathcal{L} f(g) = \lim_{s \rightarrow 0} \frac{P_s(P_t f(g)) - P_t
f(g)}{s}
$$

We now to restrict ourselves to the case of $G${\bf -invariant
L\'evy processes}, which have been recently studied in \cite{LIAO}:
\begin{df} A L\'evy process, $g_t$, is said to be {\bf
$G$-invariant} if its distribution $u_t$ is $G$-invariant.
\end{df}

The $G$-invariance ensures a sufficiently ``nice'' expression of the
transition density in terms of characters of $G$.  Let $\psi_\lambda
= \chi_\lambda / d_\lambda$ be the normalised character.  We now
have the following:
\begin{theorem}\label{371} (\cite{LIAO} Thm. 4)  Let $G$ be a compact connected Lie group and let $g_t$ be a $G$-invariant, non-degenerate L\'evy process in $G$.  Then
\begin{list}{}{}
\item (i)  For $t > 0$, the distribution $u_t$ of $g_t$ has a density
$p_t \in L^2(G)$ given by
$$
p_t(g) = \sum_{\lambda \in \Lambda_+} d_\lambda a_\lambda (t)
\chi_\lambda (g) \phantom{abcde} g \in G
$$
where $a_\lambda (t) = u_t( \bar{\psi}_\lambda ) = e^{t \mathcal{L}
(\bar{\psi}_\lambda)(e)}$, and the series converges absolutely and
uniformly for $(g,t) \in G \times \R^+$, and
$$
|a_\lambda (t)| = \exp \biggl\{ -\bigl[ \theta_\lambda + \int(1-
{\rm Re} \, \psi_\lambda) d\Pi \bigr] t \biggr\}
$$
with $\theta_\lambda = -\sum_{i,j = 1}^n a_{ij} X_i X_j
\bar{\psi}_\lambda (e) > 0$, and $\Pi$ the L\'evy measure.
\item (ii) Let
$$
\theta = \inf \biggl\{ \biggl[ \theta_\lambda + \int(1- {\rm Re} \,
\psi_\lambda) d\Pi \biggr]; \; \lambda \in \Lambda^+ \biggr\}
$$
then $\theta > 0$ for some $\lambda \in \Lambda^+$, and
$$
\| p_t - 1 \|_\infty \leq C e^{-\theta t}, \phantom{abcde}
ce^{-\theta t} \leq \| p_t - 1 \|_2 \leq C e^{-\theta t}
$$
\item (iii)  If $G$ is semisimple and the L\'evy measure has finite
first moment, then
$$
a_\lambda (t) = \exp \biggl\{ -\bigl[ \theta_\lambda + \int(1- {\rm
Re} \, \psi_\lambda) d\Pi \bigr] t \biggr\}
$$
\end{list}
\end{theorem}

In general, the wrap of $\mathcal{L}$ is difficult to determine.
Even if we consider the case where $\mathcal{L}$ is just a
differential operator, the coefficients of $\mathcal{L}$ may not be
constant (and potentially very badly behaved), and thus explicit
forms of the Duflo isomorphism may be hard to calculate. Applying
the Feynman-Ka\v{c} type transformation is complicated by the
presence of these terms. However, we are able to recover Theorem
\ref{371} (i)  -  in law  - by wrapping:
\begin{prop}\label{prop344}  Suppose $\gamma_t$ is a L\'evy process in $\g$, with distribution $h_t = \E(\gamma_t)$.  Then the distribution of the wrapped L\'evy process $\Phi(h_t)$ is given
by
$$
\Phi(h_t) (x) = \sum_{\lambda \in \Lambda^+} d_\lambda \, \hat{h_t}
(\lambda + \rho) \chi_\lambda (x), \phantom{abcde} \forall x \in T,
\, t \in \R^+.
$$
\end{prop}
\noindent {\bf Proof:}  This follows from (\ref{343}) and Lemma
\ref{PhiHat}, that $\Phi^\wedge h_t (\lambda + \rho) = \hat{h}_t
(\lambda + \rho)$ and Proposition \ref{wrapS}. $\phantom{ab}
\square$

\section{Further directions}

In the sequel, we will examine wrapping Brownian motion and heat
kernels for the cases of compact and non-compact symmetric spaces.
The wrapping formula needs some modification to hold for these more
general spaces.  This involves ``twisting" the convolution product on the tangent space
by a certain function $e$, which originates in the work of Rouvi\`ere \cite{R}.  See also \cite{D1}.\\

Ultimately, this leads us to give a concise explanation as to why
the ``sum over classical paths" (as it is known in the physics
literature) does not hold for general compact symmetric spaces
(\cite{CAM}, \cite{DOW}).\\

We will also that we have been able to extend our methods on
wrapping Brownian motion and heat kernels to some spaces where we
know the wrapping formula holds.  A nice example are the complex Lie
groups. Instead of having to deal with a maximal torus $\T^n$, as in
the case of a compact Lie group, the subgroup corresponding to the
Cartan subalgebra is $(\R^+)^n$, so instead of summing over a
lattice, we just ``bend" the heat kernel from $\g$ to $G$.  This
recovers a formula of Gangolli \cite{GAN}.\\

{\sc School of Mathematics, UNSW, Kensington 2052 NSW, Australia,

And

Group Market Risk, National Australia Bank, 24/255 George St, Sydney
2000 NSW, Australia.}\\

Email: {\tt David.G.Maher@nab.com.au}

\end{document}